\documentclass{article}
\usepackage{amssymb}
\usepackage{amsfonts}
\usepackage{graphicx}
\usepackage{amsmath}
\usepackage{booktabs}
\usepackage{multirow}
\usepackage{mathrsfs}

\newtheorem{proposition}{Proposition}
\newtheorem{remark}{Remark}

\newenvironment{proof}[1][Proof]{\textbf{#1.} }{\ \rule{0.5em}{0.5em}}

\begin{document}
\date{}
\title{A Gaussian method for the operator square root}
\author{Eleonora Denich\thanks{%
Dipartimento di Matematica e Geoscienze, Universit\`{a} di Trieste,
eleonora.denich@phd.units.it} \and Paolo Novati \thanks{%
Dipartimento di Matematica e Geoscienze, Universit\`{a} di Trieste,
novati@units.it}}
\maketitle

\begin{abstract}
We consider the approximation of the inverse square root of regularly accretive operators in Hilbert spaces. The approximation is of rational type and comes from the use of the Gauss-Legendre rule applied to a special integral formulation of the problem. We derive sharp error estimates, based on the use of the numerical range, and provide some numerical experiments.
For practical purposes, the finite dimensional case is also considered. In this setting, the convergence is shown to be of exponential type.
\end{abstract}



\section{Introduction}

Let ${\mathcal{H}}$ be a generic Hilbert space with scalar product denoted by $\left\langle \cdot ,\cdot \right\rangle $ and corresponding norm $\left\Vert u\right\Vert =\left\langle u,u\right\rangle ^{1/2}$, $u\in {\mathcal{H}}$. 
Given a linear operator ${\mathcal{L}}$ acting in ${\mathcal{H}}$, in this work we are interested in the numerical approximation of ${\mathcal{L}}^{-1/2}$, where ${\mathcal{L}}$ is assumed to be regularly accretive, that is, associated with a regular sesquilinear form (see \cite{Kato} for a background). It is known that such operators are unbounded and satisfy
\begin{equation} \label{accretive}
| \Im\left\langle {\mathcal{L}}u,u\right\rangle |
\leq \eta \Re \left\langle {\mathcal{L}}u,u\right\rangle ,
\end{equation}
for some $\eta \geq 0$, where the symbols $\Im$ and $\Re$ indicate the imaginary and the real part respectively. 
Let 
\begin{equation} \label{settore}
\Sigma_{\beta,a}=\left\{ z\in \mathbb{C} : \left\vert \arg (z-a)\right\vert \leq \beta \pi, \; a \geq 0, \; \beta < \frac{1}{2} \right\},
\end{equation}
be the sector symmetric with respect to the real axis with vertex in $a$ and semiangle $\beta \pi$.
Denoting by $F({\mathcal{L}})$ the numerical range of ${\mathcal{L}}$, that is,
\begin{equation*}
F(\mathcal{L}) = \left\lbrace z \in \mathbb{C} : z = \frac{\langle u, \mathcal{L} u \rangle}{\langle u,u \rangle}, u \in \mathcal{H}, u \neq 0 \right\rbrace ,
\end{equation*}
referring to (\ref{accretive}) it is also known that (\cite[Th. 2.2]{Kato})
\begin{equation*}
F({\mathcal{L}})\subseteq \Sigma_{\arctan (\eta)/\pi,0 }. 
\end{equation*}
In this setting, for $\alpha \in (0,1)$ the fractional power is defined by (see \cite{Bala})
\begin{equation} \label{fractional_power}
{\mathcal{L}}^{-\alpha }v=\frac{\sin (\alpha \pi) }{\pi }\int\nolimits_{0}^{\infty }\varsigma^{-\alpha }(\varsigma{\mathcal{I}}+\mathcal{L})^{-1}vd\varsigma,\quad v\in {\mathcal{H}},
\end{equation}
where ${\mathcal{I}}$ is the identity operator in ${\mathcal{H}}$. 

Starting from this representation with $\alpha =1/2$, in this work we consider some changes of variable that, in the scalar case, lead to the formula
\begin{equation} \label{caso_scalare}
\lambda^{-\frac{1}{2}} = \frac{4 \sqrt{\tau}}{\pi} \left( \int_{-1}^1 \frac{1}{4 \tau+\lambda (t+1)^2} dt + \int_{-1}^1 \frac{1}{\tau (t+1)^2+4 \lambda} dt \right),
\end{equation}
where $\tau >0$ is a parameter that allows to balance the contribute of the two integrals (Section \ref{Section4}). For the approximation of (\ref{caso_scalare}) we employ the Gauss-Legendre rule. 
Working with an $n$-point formula for both integrals we implicity construct a rational form of type
\begin{equation*}
\mathcal{R}_{2n-1,2n} (\lambda) = \frac{p_{2n-1}(\lambda)}{q_{2n}(\lambda)}, \quad p_{2n-1} \in \Pi_{2n-1}, \; q_{2n} \in \Pi_{2n},
\end{equation*}
such that 
\begin{equation*}
\mathcal{R}_{2n-1,2n} (\lambda) \cong \lambda^{-\frac{1}{2}}.
\end{equation*}
Assuming that $F(\mathcal{L}) \subset \Sigma_{\beta,a}$, $a>0$, we define $\tau$ (depending on $n$) such that the rational approximation is still reliable for $\mathcal{L}^{-\frac{1}{2}}$.

Similar approaches, based on quadrature rules arising from the
Dunford-Taylor integral representation of $\lambda ^{-\alpha }$, have been considered for instance in \cite
{ABDN,AN0,AcetoNovati,ANI,Bonito2015,Bonito2019}.
Other methods that rely on the best uniform rational approximations of functions closely related to $\lambda ^{-\alpha }$ have been treated in \cite{HLMMV, HLMMP, HML, HM}.
Methods based on the parabolic reformulation of fractional diffusion equations have been analyzed in \cite{V15,V18,V19}.
As pointed out in \cite{H} they can still be interpreted as rational approximations of $\lambda^{-\alpha}$.
We also quote here \cite{survey} for a very recent survey.
We remark that, except for \cite{Bonito2019}, in all these papers the basic assumption has been to work with a self-adjoint operator. 

As for the method considered in this paper, we are able to show that in the operator norm the convergence rate is of type 
\begin{equation*}
c\frac{\left( \ln n\right) ^{2}}{n^{4}}
\end{equation*}
where $n$ is the number of quadrature points and $c$ is a constant depending on the angle of the sector containing $F({\mathcal{L}})$. We do not claim that the proposed method is the fastest since, by using an exponential transform in (\ref{fractional_power}) and then the trapezoidal rule, it is possible to achieve an exponential decay for the error as shown in \cite{Bonito2019}. Nevertheless the method has some potentials. Indeed, the rate of convergence is independent of the angle of the sector containing $F({\mathcal{L}})$, the initial convergence is very fast because of the factor $1/n^{4}$, and,
in addition, we have been able to derive a sharp error estimate that allows to select a priori the number of quadrature points to achieve a certain accuracy.

For practical purposes we have also analyzed the behavior of the method in finite dimension, that is, in the case of bounded sectorial operators ${\mathcal{L}}_{N}$, showing that the decay of the error is of type
\begin{equation*}
\frac{c_{1}}{\left\Vert {\mathcal{L}}_{N}\right\Vert ^{1/4}}\exp \left(
-c_{2}\frac{n}{\left\Vert {\mathcal{L}}_{N}\right\Vert ^{1/8}}\right) 
\end{equation*}
where $c_{1}$ and $c_{2}$ are constants depending on the angle of the sector containing $F({\mathcal{L}}_{N})$, and $\left\Vert {\mathcal{\cdot }}\right\Vert $ is the spectral norm. It is interesting to observe that in the case of ${\mathcal{L}}_{N}$ symmetric and positive definite the above formula can be rewritten replacing $\left\Vert {\mathcal{L}}_{N}\right\Vert $ with the specral condition number $\kappa ({\mathcal{L}}_{N})$, resulting in
a clear improvement with respect for instance to the Gauss-Jacobi approach \cite{AcetoNovati}, where the above formula still holds, but with $\|{\mathcal{L}}_{N}\| $ replaced by $\|{\mathcal{L}}_{N}\|^{2}$.

The paper is organized as follow. In Section \ref{Section2} we derive the integral representation (\ref{caso_scalare}) and approximate it by using the Gauss-Legendre rule. In Section \ref{Section3} we develop the error analysis in the scalar case by studying the poles of the integrand functions in the complex plane. In Section \ref{Section4} we generalize the analysis for regularly accretive operators, by considering the error behavior on the boundaries of the sector containing the numerical range of the operator. Finally, in Section \ref{Section5}, working with bounded operators, we improve the error estimates previously obtained.

\section{The Gauss-Legendre approach} \label{Section2}

Starting from the integral representation (\ref{fractional_power}), we consider the change of variable $\varsigma=x^{-2}$ (see \cite{Bonito2015}), that leads to
\begin{equation*}
\lambda^{-\frac{1}{2}} = \frac{2}{\pi} \int_0^{\infty} \frac{1}{1+x^2 \lambda} dx, \quad \lambda \in \mathbb{C} \setminus (-\infty,0].
\end{equation*}
Then, we split the above integral as follows
\begin{equation*} 
\lambda^{-\frac{1}{2}} = \frac{2}{\pi} \left( \int_{\frac{1}{\sqrt{\tau}}}^{\infty}  \frac{1}{1+x^2 \lambda} dx + \int^{\frac{1}{\sqrt{\tau}}}_0  \frac{1}{1+x^2 \lambda} dx \right), 
\end{equation*}
where $\tau \geq 1$ is a certain parameter whose meaning will be explained later.
By using the changes of variable 
\begin{equation*}
x=\frac{1}{\sqrt{\tau} y} \quad {\rm and} \quad x=\frac{y}{\sqrt{\tau} },
\end{equation*}
for the first and the second integral respectively, we have that
\begin{equation} \label{int_spezzato2}
\lambda^{-\frac{1}{2}} = \frac{2 \sqrt{\tau}}{\pi} \left( \int_0^1 \frac{1}{\tau+\lambda y^2} dy + \int_0^1 \frac{1}{\tau y^2+\lambda} dy \right) .
\end{equation}
Finally, for both integrals in (\ref{int_spezzato2}) we apply the change of variable 
\begin{equation*}
y=\frac{t+1}{2}
\end{equation*}
and obtain the integral representation (\ref{caso_scalare}), that is,
\begin{equation} \label{AA}
\lambda^{-\frac{1}{2}}=\frac{4 \sqrt{\tau}}{\pi} \left( I^{(1)}(\lambda) + I^{(2)}(\lambda) \right),
\end{equation}
where
\begin{equation} \label{integraleA1}
I^{(1)}(\lambda):= \int_{-1}^1 \frac{1}{4 \tau+\lambda (t+1)^2} dt, \quad I^{(2)}(\lambda):=  \int_{-1}^1 \frac{1}{\tau (t+1)^2+4 \lambda} dt.
\end{equation}
Using the $n$-point Gauss-Legendre quadrature rule, the formula (\ref{AA}) is approximated as
\begin{equation} \label{A}
\lambda^{-\frac{1}{2}}\cong \frac{4 \sqrt{\tau}}{\pi} \left( I_n^{(1)}(\lambda) + I_n^{(2)}(\lambda) \right),
\end{equation}
where 
\begin{equation} \label{sommatorie1}
I_n^{(1)}(\lambda):=\sum_{j=1}^n \omega_j \left( 4 \tau+\lambda (t_j+1)^2 \right)^{-1}, \quad I_n^{(2)}(\lambda) := \sum_{j=1}^n \omega_j \left( \tau (t_j+1)^2 +4 \lambda \right)^{-1},
\end{equation}
in which $t_j$, $\omega_j$, $j=1,\ldots,n$, are respectively the nodes and the weights of the Gaussian rule.
As mentioned in the introduction, we observe that (\ref{A}) represents a rational approximation $\mathcal{R}_{2n-1,2n}(\lambda)$ of $\lambda^{-\frac{1}{2}}$, where
\begin{equation*}
\mathcal{R}_{2n-1,2n}(\lambda)=\frac{p_{2n-1}(\lambda)}{q_{2n}(\lambda)}, \quad p_{2n-1} \in \Pi_{2n-1}, \; q_{2n} \in \Pi_{2n}.
\end{equation*}

\section{General error analysis} \label{Section3}

Let us consider the transform
\begin{equation*}
\psi(w) = \frac{1}{2} \left( w+\frac{1}{w} \right), \quad w \in \mathbb{C}, \; |w|>1,
\end{equation*}
that conformally maps the exterior of the unit circle onto the exterior of the interval $[-1,1]$. This map is usually called Joukowsky transform. The image of the circle $|w| = s$, that is,
\begin{equation*}
\Psi_s = \left\lbrace z \in \mathbb{C} \colon z=\psi \left(se^{i\theta}\right), \, \theta \in [0, 2 \pi] \right\rbrace,
\end{equation*}
is an ellipse of the complex plane with foci in $\pm 1$.
We denote by $\mathcal{E} = \left\{ \Psi_s \, | \, s>1 \right\}$ the family of all these ellipses.

\begin{figure}
\begin{center}
\includegraphics[scale=0.3]{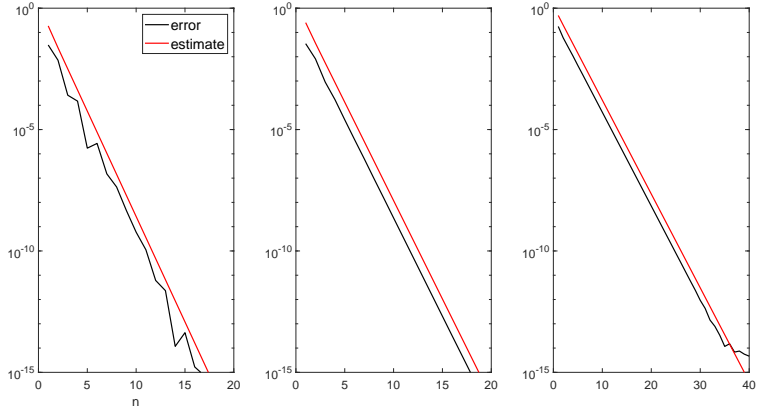}
\end{center}
\caption{Error and error estimate (\ref{moduloE_n}) for $\lambda=10$, $\lambda=5+5i$, $\lambda=-5+10i$, from left to right, and $\tau=2$.}
\label{Figura1}
\end{figure}

Let $f$ be a generic function analytic in an open set containing $[-1,1]$. Let moreover 
\begin{equation*}
I(f) = \int_{-1}^1 f(t) dt,
\end{equation*}
and $I_n(f)$ its $n$-point Gauss-Legendre approximation. Following the analysis given in \cite{Barrett}, assume that $\Psi_s \in \mathcal{E}$ is such that $f$ is analytic in the interior of $\Psi_s$, except for a pair of simple poles $t_0$ and its conjugate $\bar{t}_0$. Then, 
\begin{equation} \label{B}
I(f)-I_n(f) \cong -4 \pi \Re \left\lbrace r \left( t_0 + \sqrt{t_0^2-1} \right)^{-2n} \right\rbrace,
\end{equation} 
where $r$ is the residue of $f$ at $t_0$ and the root has to be choosen such that 
\begin{equation*}
S:=\Big| t_0 + \sqrt{t_0^2-1} \Big| >1.
\end{equation*}
It is important to observe that $S$ is just the radius of the circle centered at $0$ that, through the Joukowsky transform, is mapped onto the ellipse $\Psi_S \in \mathcal{E}$ passing through $t_0$ and $\bar{t}_0$. Clearly, by (\ref{B}) the rate of convergence grows with $S$, that roughy speaking, handle the distance of the poles from the interval $[-1,1]$.
Having at disposal the above general result, we can estimate the error of the approximation (\ref{A}) by studying separately the poles of our integrand functions  (cf. (\ref{integraleA1}))
\begin{equation} \label{funzioni}
f^{(1)}(t) = \frac{1}{4 \tau +\lambda (t+1)^2}
\quad {\rm and} \quad
f^{(2)}(t) = \frac{1}{\tau (t+1)^2+4 \lambda}.
\end{equation}
As for the function $f^{(1)}$ it is easy to see that the poles are given by 
\begin{equation} \label{poli1}
t_{0,1}^{(1)} = \pm 2 \left( \frac{\tau}{\lambda} \right)^{\frac{1}{2}} i -1, 
\end{equation}
so that $t_0^{(1)}= \overline{t_1^{(1)}}$.
Similarly, for the function $f^{(2)}$ we have 
\begin{equation} \label{poli2}
t_{0,1}^{(2)} = \pm 2 \left( \frac{\lambda}{\tau} \right)^{\frac{1}{2}} i -1, 
\end{equation}
and therefore $t_0^{(2)}= \overline{t_1^{(2)}}$.
By using (\ref{B}), and defining $e_n^{(i)}(\lambda) =  I^{(i)}(\lambda)- I_n^{(i)}(\lambda) $, $i=1,2$ (cf. (\ref{AA}) and (\ref{A})), we obtain
\begin{equation*}
e_n^{(i)}(\lambda) \cong -4 \pi \Re \left\lbrace r^{(i)} \left( t_0^{(i)} + \sqrt{\left(t_0^{(i)}\right)^2-1} \right)^{-2n} \right\rbrace.
\end{equation*}
Since 
\begin{equation*}
r^{(i)} = \lim_{t \rightarrow t_0^{(i)}} \left( t-t_0^{(i)} \right) f^{(i)}(t), \quad i=1,2,
\end{equation*}
by (\ref{funzioni}), (\ref{poli1}), (\ref{poli2}) we find
\begin{equation*}
r^{(1)} = \frac{1}{4 \sqrt{\tau}} \lambda ^{-\frac{1}{2}} e^{i \frac{\pi}{2}} \quad {\rm and} \quad r^{(2)} = \frac{1}{4 \sqrt{\tau}} \lambda ^{-\frac{1}{2}} e^{-i \frac{\pi}{2}},
\end{equation*}
and therefore
\begin{align}
|e_n^{(i)}(\lambda)| &\cong 4 \pi |r^{(i)}| {S^{(i)}}^{-2n} \notag \\
&= \frac{\pi}{\sqrt{\tau}} |\lambda|^{-\frac{1}{2}} {S^{(i)}}^{-2n} =: \Phi^{(i)} (\tau,\lambda), \label{PHI}
\end{align}
where 
\begin{equation} \label{k_i}
S^{(i)}=\Bigg|t_0^{(i)} + \sqrt{\left(t_0^{(i)}\right)^2-1}\Bigg|, \quad i=1,2.
\end{equation}
As for the total error
\begin{equation} \label{E_n}
E_n(\lambda) = \frac{4 \sqrt{\tau}}{\pi} \left( e_n^{(1)}(\lambda) + e_n^{(2)}(\lambda) \right),
\end{equation}
(cf. (\ref{AA}) and (\ref{A})), we then have the estimate
\begin{equation} \label{moduloE_n}
|E_n(\lambda)|  \cong  4 |\lambda|^{-\frac{1}{2}} \left( {S^{(1)}} ^{-2n}+  {S^{(2)}} ^{-2n} \right).
\end{equation}
In Figure \ref{Figura1} we show the accuracy of the above formula for some values of $\lambda$ and $\tau =2$.

\section{Estimates for operators} \label{Section4}

\begin{figure} 
\begin{center}
\includegraphics[scale=0.3]{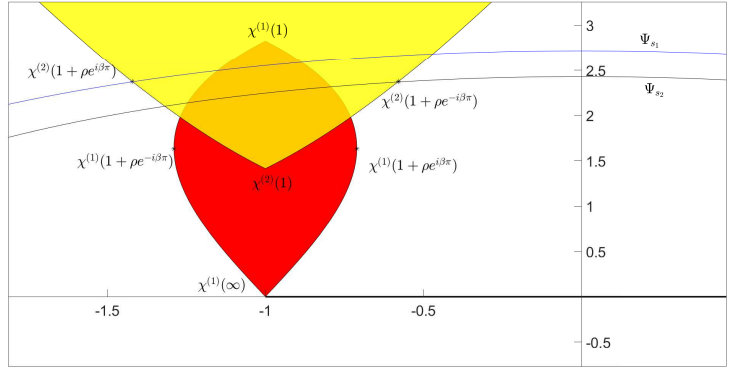}
\end{center}
\caption{The regions of the poles (with positive imaginary part) of the functions $f^{(1)}(t)$ (red) and $f^{(2)}(t)$ (yellow), for $\beta = \frac{1}{6}$ and $\tau=2$, and the ellipses $\Psi_{s_1}$ and $\Psi_{s_2}$.}
\label{Figura2}
\end{figure}

Let $\mathcal{H}$ be a generic Hilbert space and let $\mathcal{L} \colon \mathcal{H}\rightarrow \mathcal{H}$ be a regularly accretive operator such that (see (\ref{settore}))
\begin{equation*}
F(\mathcal{L}) \subseteq \Sigma_{\beta,1}.
\end{equation*}
It is known that for any function $f$ analytic in $F(\mathcal{L})$ it holds
\begin{equation*}
\| f(\mathcal{L}) \|_{\mathcal{H}\rightarrow \mathcal{H}} \leq K \max_{\lambda \in F(\mathcal{L})} | f(\lambda) |,
\end{equation*}
where $2 \leq K \leq 1+\sqrt{2}$ is the absolute constant studied in \cite{Crouzeix}. 
We remark that if $\mathcal{L}$ is self-adjoint, then $\beta=0$ and $K=1$.
As consequence, since the poles of the approximation $\mathcal{R}_{2n-1,2n} (\lambda)$ are all in $\mathbb{R}^-$ (cf. (\ref{A}) and (\ref{sommatorie1})), we can consider the bound 
\begin{equation*} 
\|E_n(\mathcal{L})\|_{\mathcal{H}\rightarrow \mathcal{H}}  = \| \mathcal{L}^{-\alpha} - \mathcal{R}_{2n-1,2n}(\mathcal{L}) \|_{\mathcal{H}\rightarrow \mathcal{H}} \leq K \max_{\lambda \in \Sigma_{\beta,1}} | \lambda^{-\alpha}-\mathcal{R}_{2n-1,2n}(\lambda) |.
\end{equation*}
When studying the behavior of the method applied to $\lambda^{-\frac{1}{2}}$, for $\lambda \in \Sigma_{\beta,1}$, that is $\lambda=1+\rho e^{i \theta \pi}$, $\rho \geq 0$, $| \theta | \leq \beta$, it is rather evident (and will be confirmed by the following analysis) that, for a fixed $\rho$, moving $\theta$ from $0$ to $\beta$ causes a progressive slow down. 
Further, since $\lambda^{-\alpha}- \mathcal{R}_{2n-1,2n}(\lambda)$ is analytic in $\Sigma_{\beta,1}$, by the maximum modulus principle it is sufficient to consider the scalar error on 
\begin{equation*}
\Gamma_{\beta}=\Gamma_{\beta}^+ \cup \Gamma_{\beta}^-,
\end{equation*}
where
\begin{equation*}
\Gamma_{\beta}^+ = \left\lbrace z \in \mathbb{C} \mid z=1+\rho e^{i \beta \pi}, \rho \geq 0 \right\rbrace,
\end{equation*}
and 
\begin{equation*}
\Gamma_{\beta}^- = \left\lbrace z \in \mathbb{C} \mid z=1+\rho e^{-i \beta \pi}, \rho \geq 0 \right\rbrace.
\end{equation*}
Therefore, by using (\ref{E_n}) we have that
\begin{align*}
\max_{\lambda \in \Sigma_{\beta,1}} | \lambda^{-\alpha}-\mathcal{R}_{2n-1,2n}(\lambda) | &\leq \max_{\lambda \in \Gamma_{\beta}} | \lambda^{-\alpha}-\mathcal{R}_{2n-1,2n}(\lambda) | \\
&\leq \frac{4 \sqrt{\tau}}{\pi} \left( \max_{\lambda \in \Gamma_{\beta}} |e_n^{(1)}(\lambda)| + \max_{\lambda \in \Gamma_{\beta}} |e_n^{(2)}(\lambda)| \right).
\end{align*}
At this point, for $\lambda \in \Gamma_{\beta}$, we consider (cf. (\ref{poli1}) and (\ref{poli2})) the maps
\begin{equation*}
\chi^{(1)} (\lambda) := 2 \left(\frac{\tau}{\lambda}\right)^{\frac{1}{2}}i-1
\end{equation*}
and
\begin{equation*}
\chi^{(2)} (\lambda) := 2 \left(\frac{\lambda}{\tau}\right)^{\frac{1}{2}}i-1,
\end{equation*}
which define the boundaries of the regions of the poles (with positive imaginary part) of the functions $f^{(1)}(t)$ and $f^{(2)}(t)$, respectively. 
These regions are plotted in Figure \ref{Figura2}. The symmetry of $\Gamma_{\beta}$ with respect to the real axis leads to the symmetry of $\chi^{(i)}(\Gamma_{\beta})$ with respect to the line $\Re (z) =-1$.
Indeed, for each fixed $\rho>0$, the two points $\chi^{(i)}(1+\rho e^{\pm i \beta \pi})$, $i=1,2$, are symmetric with respect to the line $\Re (z) =-1$. 
The ellipse passing through the one with the real part greater than $-1$ is the smallest, so that, as already observed, it leads to the worst case in terms of rate of convergence (cf. (\ref{B})). Just for clarity, in Figure \ref{Figura2} we also plot the two ellipses $\Psi_{s_1}, \, \Psi_{s_2} \in \mathcal{E}$ passing through $\chi^{(2)}(1+\rho e^{\pm i \beta \pi})$, where
\begin{align*}
s_1 &= \left\vert \chi^{(2)}(1+\rho e^{ i \beta \pi})+\sqrt{\left(\chi^{(2)}(1+\rho e^{i \beta \pi})\right)^2-1} \right\vert, \\
s_2 &= \left\vert \chi^{(2)}(1+\rho e^{- i \beta \pi})+\sqrt{\left(\chi^{(2)}(1+\rho e^{-i \beta \pi})\right)^2-1} \right\vert.
\end{align*}
As consequence, since $\Re \left( \chi^{(1)}(\lambda) \right) \geq -1 $, for $\lambda \in \Gamma_{\beta}^+$, and $\Re \left( \chi^{(2)}(\lambda) \right) \geq -1 $, for $\lambda \in \Gamma_{\beta}^-$, we have that
\begin{equation*}
\max_{\lambda \in \Gamma_{\beta}} |e_n^{(1)}(\lambda)| = \max_{\lambda \in \Gamma_{\beta}^+} |e_n^{(1)}(\lambda)|
\end{equation*}
and
\begin{equation*}
\max_{\lambda \in \Gamma_{\beta}} |e_n^{(2)}(\lambda)| = \max_{\lambda \in \Gamma_{\beta}^-} |e_n^{(2)}(\lambda)|.
\end{equation*}
Therefore, we can finally write 
\begin{align}
\max_{\lambda \in \Sigma_{\beta}} | \lambda^{-\alpha}-R_{2n-1,2n}(\lambda) | &\leq \frac{4 \sqrt{\tau}}{\pi} \max_{\rho \geq 0}\left( |e_n^{(1)}(1+\rho e^{i\beta\pi})| + |e_n^{(2)}(1+\rho e^{-i\beta\pi})| \right) \notag \\
&\cong \frac{4 \sqrt{\tau}}{\pi} \max_{\rho \geq 0} \left( \Phi^{(1)}(\tau,1+\rho e^{i\beta\pi}) + \Phi^{(2)}(\tau,1+\rho e^{-i\beta\pi}) \right), \label{maxmax}
\end{align}
where we have used the relation (\ref{PHI}).

\subsection{Error behavior on $\Gamma_{\beta}$}

Experimentally, working with operators with large spectrum, one observes that, in general, the parameter $\tau$ must be choosen quite large to achieve a good rate of convergence. For this reason from now we assume $\tau \gg 1$. By considering the functions
\begin{equation*}
\Phi^{(1)}\left(\tau, 1+\rho e^{i \beta \pi} \right) \quad {\rm and} \quad \Phi^{(2)}\left(\tau, 1+\rho e^{-i \beta \pi} \right),
\end{equation*}
for $\rho \geq 0$ (cf. (\ref{maxmax}),(\ref{PHI})), it can be observed that, with respect to $\rho$, $\Phi^{(1)}\left(\tau, 1+\rho e^{i \beta \pi} \right)$ initially grows, reach a maximum at a certain $\rho^{\star} \gg \tau$, and then goes to $0$ as $\rho \rightarrow +\infty$. On the other side, $\Phi^{(2)}\left(\tau, 1+\rho e^{-i \beta \pi} \right)$ may show two kinds of behavior, depending on the angle $\beta \pi$. In particular, for $\beta \leq \beta^{\star}$ ($\beta^{\star} \cong \frac{1}{4}$), $\Phi^{(2)}\left(\tau, 1+\rho e^{-i \beta \pi} \right)$ is monotone decreasing, whereas for $\beta > \beta^{\star}$ it initially grows, reach a maximun at a certain $\bar{\rho} \ll \tau$, and then is monotone decreasing. 
We refere to the Appendix for the details concerning $\beta^{\star}$. 
Similarly to the analysis given in \cite{AcetoNovati}, the idea is then to define $\tau$ in order that
\begin{equation} \label{condizionePHI}
\Phi^{(1)}\left(\tau, 1+\rho^{\star} e^{i \beta \pi} \right) = \Phi^{(2)}(\tau, 1+\bar{\rho}e^{-i \beta \pi}),
\end{equation}
where we should set $\bar{\rho}=0$ for $\beta \leq \beta^{\star}$.
\begin{remark} \label{Remark1}
The reason why we impose (\ref{condizionePHI}) is that numerically one observes that (cf. (\ref{maxmax}))
\begin{equation*}
\Phi^{(1)}\left(\tau, 1+\rho e^{i \beta \pi} \right)+\Phi^{(2)}\left(\tau, 1+\rho e^{-i \beta \pi} \right) \cong \max_{\rho \geq 0} \left( \Phi^{(1)}\left(\tau, 1+\rho e^{i \beta \pi} \right),\Phi^{(2)}\left(\tau, 1+\rho e^{-i \beta \pi} \right)\right).
\end{equation*}
\end{remark}

\subsubsection{Approximation of $\rho^{\star}$}

In what follows we use the symbol $\sim$ to relate functions asymptotically equal in the usual sense. Since $\rho^{\star} \gg 1$, in order to study the function $\Phi^{(1)}\left(\tau, 1+\rho e^{i \beta \pi} \right)$, we first consider the approximation
\begin{align}
t_0^{(1)} &= 2 \left( \frac{\tau}{1+\rho e^{i \beta \pi}} \right) ^{\frac{1}{2}} i -1 \notag \\
&\sim 2 \left( \frac{\tau}{\rho} \right) ^{\frac{1}{2}} e^{i \frac{\pi}{2} (1- \beta)} -1, \quad \rho \rightarrow \infty. \label{t_0 approx}
\end{align} 
As consequence, for the term $S^{(1)}$ (see (\ref{PHI}) and (\ref{k_i})) we have the following result.
\begin{proposition}
It holds
\begin{equation} \label{k_1 approx}
S^{(1)} \sim 1+ \sqrt{2} C_{\beta} \left( \frac{\tau}{\rho} \right) ^{\frac{1}{4}}, \quad \rho \rightarrow \infty,
\end{equation}
where 
\begin{equation} \label{C_beta}
C_{\beta} = \sqrt{2} \cos \left(\frac{\pi}{4} (\beta +1) \right).
\end{equation}
\end{proposition}
\begin{proof}
Assuming that $\rho \gg \tau$, by (\ref{k_i}), (\ref{t_0 approx}) and using the first order approximation 
\begin{equation} \label{radice}
\sqrt{1+x} = 1+ \frac{x}{2}+\mathcal{O}\left( x^2 \right), \quad {\rm for} \quad x \rightarrow 0,
\end{equation}
we can write
\begin{align*}
S^{(1)} & \sim \Bigg | 2 \left( \frac{\tau}{\rho} \right) ^{\frac{1}{2}} e^{i \frac{\pi}{2} (1- \beta)} -1 +\sqrt{4 \left( \frac{\tau}{\rho} \right) e^{i \pi (1- \beta)} -4 \left( \frac{\tau}{\rho} \right) ^{\frac{1}{2}} e^{i \frac{\pi}{2} (1- \beta)}} \Bigg| \\
&\sim \Bigg| -1+2i \left( \frac{\tau}{\rho} \right)^\frac{1}{4} e^{i \frac{\pi}{4} (1- \beta)} \Bigg| \\
&= \sqrt{\left[-1+\sqrt{2} \left( \frac{\tau}{\rho} \right)^{\frac{1}{4}} \left(\sin \left(\frac{\beta}{4}\pi \right)-\cos \left(\frac{\beta}{4}\pi \right) \right) \right]^2 +2 \left( \frac{\tau}{\rho} \right)^{\frac{1}{2}} \left(\sin \left( \frac{\beta}{4}\pi \right) +\cos \left( \frac{\beta}{4} \pi \right) \right)^2} \\
&\sim \sqrt{1+2\sqrt{2} \left( \frac{\tau}{\rho} \right)^{\frac{1}{4}} \left( \cos \left(\frac{\beta}{4}\pi \right) -\sin \left( \frac{\beta}{4} \pi \right) \right)},
\end{align*}
that leads to (\ref{k_1 approx}) because
\begin{equation*}
\cos \left(\frac{\beta}{4}\pi \right) -\sin \left( \frac{\beta}{4} \pi \right) = \sqrt{2} \cos \left(\frac{\pi}{4} (\beta +1) \right).
\end{equation*}
\end{proof}

Note that $0.54 \cong \sqrt{1-\frac{\sqrt{2}}{2}} < C_{\beta} \leq 1$ for $0 \leq \beta < \frac{1}{2}$ ($C_{\beta}=1$ for $\beta=0$).
At this point, we look for the local maximun of the approximation (see (\ref{PHI}))
\begin{equation} \label{g1}
\Phi^{(1)}(\tau,1+\rho e^{i\beta \pi}) \sim \frac{\pi}{\sqrt{\tau}} \rho^{-\frac{1}{2}}\left[1+\sqrt{2} C_{\beta} \left(\frac{\tau}{\rho} \right)^{\frac{1}{4} } \right]^{-2n}=:g^{(1)} (\tau,\rho),
\end{equation}
where we have also used
\begin{equation*}
\Big|\left(1+\rho e^{i \beta \pi} \right) ^{-\frac{1}{2}} \Big |\sim \rho^{-\frac{1}{2}}.
\end{equation*}
By solving
\begin{equation*}
\frac{d}{d \rho} \rho^{-\frac{1}{2}} \left[1+\sqrt{2} C_{\beta} \left(\frac{\tau}{\rho} \right)^{\frac{1}{4} } \right]^{-2n}=0,
\end{equation*}
after some computations we obtain
\begin{equation} \label{rho_star}
\hat{\rho} = 4 C_{\beta}^4 \tau (n-1)^4 \cong \rho^{\star}.
\end{equation}

\subsubsection{Approximation of $\bar{\rho}$}

\begin{figure}
\begin{center}
\includegraphics[scale=0.3]{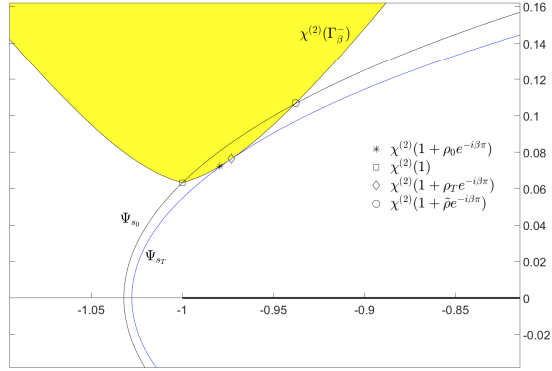}
\end{center}
\caption{Geometrical interpretation of the local maximum of $\Phi^{(2)}(\tau,1+\rho e^{-i\beta \pi})$ and its approximation by $\rho_0$, for $\beta = \frac{5}{12}$ and $\tau=1000$.}
\label{Figura6}
\end{figure}

For $\beta > \beta^{\star}$, the maximun of $\Phi^{(2)}(\tau,1+\rho e^{-i\beta \pi})$ can be approximated by considering the ellipse $\Psi_{s_T} \in \mathcal{E}$ tangent to the curve $\chi^{(2)}(\Gamma_{\beta}^-)$. 
This because the corresponding $s_T$ represents the slowest convergence rate by (\ref{PHI}). 
Let $\rho_T$ be such that $\chi^{(2)}\left(1+\rho_T e^{-i \beta \pi} \right)$ is the tangent point (see Figure \ref{Figura6}).
Since the computation of $\rho_T$ involves the solution of a $4$-th degree equation, we consider an approximation arising from geometrical evidences. We first look for the ellipse $\Psi_{s_0} \in \mathcal{E}$ ($s_0>s_T$) passing through the point $ \chi^{(2)}(1)= \frac{2}{\sqrt{\tau}}i-1$. Hence, we need to solve with respect to $s$ and $\varphi$
\begin{equation*}
\frac{1}{2} \left( s e^{i\varphi}+\frac{1}{s e ^{i\varphi}} \right) = \frac{2}{\sqrt{\tau}} i -1,
\end{equation*} 
or equivalently
\begin{equation*}
\begin{cases}
\frac{1}{2} \left( s \cos \varphi+\frac{1}{s}\cos \varphi \right) = -1 \\
\frac{1}{2} \left( s \sin \varphi-\frac{1}{s}\sin \varphi \right) = \frac{2}{\sqrt{\tau}}
\end{cases}.
\end{equation*}
After some computations, we find that the solution $(s_0,\varphi_0)$ is such that
\begin{align}
\sin \varphi_0 &= \sqrt{\frac{2}{\tau}} \sqrt{\sqrt{1+\tau}-1},  \label{seno_psi_0} \\
\cos\varphi_0 &= - \sqrt{1-\frac{2}{\tau} (\sqrt{1+\tau}-1)}, \notag
\end{align} 
with $\frac{\pi}{2}< \varphi_0 < \pi$, and
\begin{equation} \label{22_bis}
s_0 = \frac{\sqrt{2}}{\sqrt{\sqrt{1+\tau}}-1} +\frac{1}{\sqrt{1-\frac{2}{\tau}(\sqrt{1+\tau}-1)}}.
\end{equation}
The idea is to approximate $\bar{\rho}$ by looking for the other intersection between $\Psi_{s_0}$ and $\chi^{(2)}(\Gamma_{\beta}^-)$ (see again Figure \ref{Figura6}).
In particular we need to solve, with respect to $\rho$, the equation
\begin{equation*}
\frac{1}{2} \left( s_0 e^{i\varphi}+\frac{1}{s_0 e ^{i\varphi}} \right) = 2 \left(\frac{1+\rho e^{-i \beta \pi}}{\tau} \right)^{\frac{1}{2}}i -1.
\end{equation*}
Setting for simplicity $a=s_0+\frac{1}{s_0}$ and $b=s_0-\frac{1}{s_0}$, the above equation leads to the system 
\begin{equation*} 
\begin{cases}
\left(1+\frac{1}{2} a \cos \varphi\right)^2 -\frac{1}{4} b^2 \sin^2 \varphi =-\frac{4}{\tau} \left(1+\rho \cos (\beta \pi) \right)\\
b \sin \varphi \left(1+\frac{1}{2} a \cos \varphi \right) = \frac{4}{\tau} \rho \sin (\beta \pi).
\end{cases}
\end{equation*}
Substituting
\begin{equation*}
\left(1+\frac{1}{2} a \cos \varphi \right)^2 =  \frac{16}{\tau^2}\frac{\rho^2 \sin^2 (\beta \pi)}{b^2 \sin^2 \varphi}
\end{equation*}
in the first equation, we obtain
\begin{equation*}
\frac{16}{\tau^2} \frac{\rho^2 \sin^2 (\beta \pi)}{b^2 \sin^2 \varphi}-\frac{1}{4} b^2 \sin^2 \varphi +\frac{4}{\tau} (1+\rho \cos (\beta \pi)) =0,
\end{equation*}
from which, after some computations, we find the solution
\begin{equation*}
\tilde{\rho} = \frac{\tau b^2 \sin^2 \varphi}{8 \sin^2 (\beta \pi)} \left[ -\cos (\beta \pi) \pm \sqrt{1-\frac{16 \sin^2 (\beta \pi)}{\tau b^2 \sin^2 \varphi}} \right].
\end{equation*}
By using (\ref{seno_psi_0}) and (\ref{22_bis}), we have that
\begin{equation*}
b^2 =\left(s_0-\frac{1}{s_0}\right)^2= \frac{8}{\sqrt{1+\tau}-1} \sim \frac{8}{\sqrt{\tau}}, \quad \tau \rightarrow \infty.
\end{equation*}
Using again (\ref{radice}), we find
\begin{equation*}
\tilde{\rho} \sim \frac{\tau^{\frac{1}{2}}\sin^2 \varphi}{\sin^2 (\beta \pi)} \left[ - \cos (\beta \pi) \pm \left(1-\frac{\sin^2 (\beta \pi)}{\tau^{\frac{1}{2}} \sin^2 \varphi} \right) \right].
\end{equation*}
Since the angle $\varphi$ is still unknown and its computation requires again the solution af a $4$-th degree equation, by taking the positive solution, we assume that (cf. (\ref{seno_psi_0}))
\begin{equation*}
\sin^2 \varphi \cong \sin^2 \varphi_0 \sim \frac{2}{\tau^{\frac{1}{2}}},
\end{equation*}
to finally obtain the rough approximation
\begin{equation} \label{rho_0}
\tilde{\rho} \cong \rho_0 := \frac{1-\cos (\beta \pi)}{1+\cos (\beta \pi)} = \tan^2 \left(\frac{\beta \pi}{2}\right).
\end{equation}
Experimentally, we observe that $0\leq \rho_0 < \tilde{\rho}$ and that $\chi^{(2)}(1+\rho_0e^{-i\beta\pi})$, for large $\tau$, is close to the tangent point independently of $\beta$. Therefore, we use $\rho_0$ as an approximation of $\bar{\rho}$.
As for the term $S^{(2)}$ in (\ref{k_i}) we have the following result.

\begin{proposition}
For $\rho = \rho_0$ and $\beta>\beta^{\star}$, we have
\begin{equation*}
S^{(2)} \sim 1+\sqrt{2} G_{\beta} \tau^{-\frac{1}{4}}, \quad \tau \rightarrow \infty,
\end{equation*}
where 
\begin{equation} \label{G_beta}
G_{\beta}=\sqrt{D_{\beta} - \sqrt{A^-_{\beta}}},
\end{equation}
with
\begin{equation} \label{C_bar_beta}
D_{\beta} = \left( 1+2\rho_0 \cos (\beta \pi) + \rho_0^2 \right)^{\frac{1}{4}}, 
\end{equation}
and
\begin{equation} \label{A-}
A^-_{\beta}= \frac{-1-\rho_0 \cos (\beta \pi) +\sqrt{1+2 \rho_0 \cos (\beta \pi) +\rho_0^2}}{2}.
\end{equation}
\end{proposition}
\begin{proof}
First of all, for $\rho =\rho_0$  we have that
\begin{align*}
t_0^{(2)} &= \frac{2}{\sqrt{\tau}} \left( 1+\rho_0 \cos (\beta \pi) -i \rho_0 \sin (\beta \pi) \right)^{\frac{1}{2}}-1 \\
&= \frac{2}{\sqrt{\tau}} \left( i \sqrt{A^+_{\beta}}+ \sqrt{A^-_{\beta}} \right)-1,
\end{align*}
where
\begin{equation*}
A^+_{\beta}= \frac{1+\rho_0 \cos (\beta \pi) +\sqrt{1+2 \rho_0 \cos (\beta \pi) +\rho_0^2}}{2},
\end{equation*}
and $A^-_{\beta}$ as in (\ref{A-}). Now, defining $D_{\beta}:= \left( 1+2\rho_0 \cos (\beta \pi) + \rho_0^2 \right)^{\frac{1}{4}}$ and using (\ref{radice}) for large $\tau$
\begin{align*}
S^{(2)} &= \Bigg | \frac{2}{\sqrt{\tau}} \left( i \sqrt{A^+_{\beta}}+ \sqrt{A^-_{\beta}} \right)-1 +\sqrt{\frac{4}{\tau}\left( A^-_{\beta}-A^+_{\beta}+2i \sqrt{A^-_{\beta} A^+_{\beta}}\right) -\frac{4}{\sqrt{\tau}} \left(\sqrt{A^-_{\beta}}+i \sqrt{A^+_{\beta}} \right)} \Bigg | \\
& \sim \Bigg | -1+\sqrt{-\frac{4}{\sqrt{\tau}}\left(\sqrt{A^-_{\beta}}+i \sqrt{A^+_{\beta}}\right)} \Bigg | \\
&= \Bigg | -1-\frac{\sqrt{2}}{\tau^{\frac{1}{4}}} \sqrt{-\sqrt{A^-_{\beta}}+D_{\beta}}+i \frac{\sqrt{2}}{\tau^{\frac{1}{4}}} \sqrt{\sqrt{A^-_{\beta}}+D_{\beta}} \Bigg | \\
&= \sqrt{1+\frac{2}{\tau^{\frac{1}{2}}}\left( -\sqrt{A^-_{\beta}}+D_{\beta} \right) +\frac{2 \sqrt{2}}{\tau^{\frac{1}{4}}}\sqrt{ -\sqrt{A^-_{\beta}}+D_{\beta} }+\frac{2}{\tau^{\frac{1}{2}}} \left( \sqrt{A^-_{\beta}}+D_{\beta} \right)} \\
&\sim \sqrt{1+\frac{2 \sqrt{2}}{\tau^{\frac{1}{4}}} \sqrt{D_{\beta}-\sqrt{A^-_{\beta}}}}.
\end{align*}
Finally, using again (\ref{radice}) and defining
\begin{equation*}
G_{\beta} := \sqrt{D_{\beta}-\sqrt{A^-_{\beta}}},
\end{equation*}
we obtain the result.
\end{proof}

We notice that $A_{\beta}^-=0$, $D_{\beta}=1$ and therefore $G_{\beta}=1$ for $\beta =0$.
By using the above proposition and
\begin{equation*}
|1+\rho_0 e^{-i \beta \pi}| = \sqrt{1+2\rho_0 \cos (\beta \pi) + \rho_0^2} = D_{\beta}^2,
\end{equation*}
(cf. (\ref{C_bar_beta})), we obtain (see (\ref{PHI}))
\begin{equation} \label{C}
\Phi^{(2)}(\tau,1+\rho_0 e^{-i \beta \pi}) \sim \frac{\pi}{\sqrt{\tau} D_{\beta}} \left(1+\sqrt{2}G_{\beta} \tau^{-\frac{1}{4}} \right)^{-2n}=:g^{(2)}(\tau,\rho_0).
\end{equation}

\subsection{The optimal value for $\tau$} \label{Section4.2}

Working with the approximations (\ref{g1})-(\ref{rho_star}) and (\ref{C}), in order to find a nearly optimal value for $\tau$, we impose the condition
\begin{equation} \label{condizione tau}
g^{(1)}(\tau, \hat{\rho}) = g^{(2)}(\tau,\rho_0).
\end{equation}
\begin{remark}
For $0\leq \beta \leq \beta^{\star}$ we should replace $g^{(2)}(\tau,\rho_0)$ by $g^{(2)}(\tau,0)$ in (\ref{condizione tau}). Anyway, since $\rho_0 \lesssim 0.17$ for $0 \leq \beta \leq \beta^{\star} \cong \frac{1}{4}$ , we still work with (\ref{condizione tau}), because this choice does not influence the results. Observe moreover that (\ref{rho_0}) luckily gives the correct value $\rho_0=0$ for $\beta=0$.
\end{remark}
The equation (\ref{condizione tau}) leads to 
\begin{equation*} 
\frac{\tau^{-\frac{1}{2}}}{2 C_{\beta}^2 (n-1)^2} \left(1+\frac{1}{n-1} \right)^{-2n} = \frac{1}{D_{\beta}} \left( 1+\sqrt{2}G_{\beta} \tau^{-\frac{1}{4}} \right)^{-2n}.
\end{equation*}
Since
\begin{equation} \label{e^-2}
\left(1+\frac{1}{n-1} \right)^{-2n} \sim e^{-2},
\end{equation}
after some computation we rewrite the above equation as
\begin{equation*}
\ln \left(\tau^{-\frac{1}{4}} \right) + \ln \left( \frac{D_{\beta}^{\frac{1}{2}}}{\sqrt{2}C_{\beta}e(n-1)} \right) = -n \ln \left(1+\sqrt{2}G_{\beta} \tau^{-\frac{1}{4}} \right).
\end{equation*}
Let us denote by $W(x)$ the Lambert-$W$ function, for which it holds
\begin{equation*}
\frac{x}{W(x)}=e^{W(x)}.
\end{equation*}
By using 
\begin{equation*}
\ln \left(1+\sqrt{2}G_{\beta} \tau^{-\frac{1}{4}} \right) \sim \sqrt{2}G_{\beta} \tau^{-\frac{1}{4}}, \quad \tau \gg 1,
\end{equation*}
we have that an approximated optimal value for $\tau$ is given by
\begin{equation} \label{tau_star}
\bar{\tau} = \frac{D_{\beta}^2}{4C_{\beta}^4e^4(n-1)^4}  \exp{\left(4W\left( H_{\beta}n(n-1) \right)\right)},
\end{equation}
where 
\begin{equation*}
H_{\beta} := \frac{2e C_{\beta} G_{\beta}}{D_{\beta}^{\frac{1}{2}}}.
\end{equation*}
Note that by (\ref{C_beta}), (\ref{G_beta}), (\ref{C_bar_beta}), (\ref{A-}), we have $H_{\beta}=2e$ for $\beta=0$.

\subsection{Asymptotic expression of the global error}

By substituting expression (\ref{tau_star}) in (\ref{rho_star}) we obtain
\begin{equation} \label{rho_star_exp}
\hat{\rho}=\frac{D_{\beta}^2}{e^4}  \exp{\left(4W\left( H_{\beta}n(n-1) \right)\right)}.
\end{equation}
Since for large $x$ (see \cite{LambertW})
\begin{equation*}
W(x) \sim \ln x - \ln ( \ln x),
\end{equation*}
we have that 
\begin{equation*}
\exp{\left(4W\left( H_{\beta}n(n-1) \right)\right)} \sim 
\left[ \frac{ H_{\beta} n^2}{ \ln \left(  H_{\beta} n^2 \right)} \right]^{4},
\end{equation*}
and thus
\begin{equation} \label{rho_star_approx}
\hat{\rho} \sim \frac{D_{\beta}^2}{e^4}\left[ \frac{ H_{\beta} n^2}{ \ln \left(  H_{\beta} n^2 \right)} \right]^{4} .
\end{equation}
Now, since $\hat{\rho}$ rapidly grows with $n$, as pointed out in Remark \ref{Remark1} we can estimate the global error as
\begin{equation} \label{Ecerchiato}
|E_n(\mathcal{L})|  \cong  4 K \left( \hat{\rho}\right)^{-\frac{1}{2}} {S^{(1)}}^{-2n}.
\end{equation}
Finally, substituting (\ref{rho_star_approx}) in the above expression  and using again the approximation (\ref{e^-2}), after some computation we find
\begin{equation} \label{bound1}
|E_n(\mathcal{L})|  \cong  4 K \left[\frac{\ln \left(  H_{\beta} n^2 \right)}{2eC_{\beta}G_{\beta}} \right]^2 n^{-4}.
\end{equation}

\subsection{Numerical example} \label{Section4.4}

\begin{figure}
\begin{center}
\includegraphics[scale=0.3]{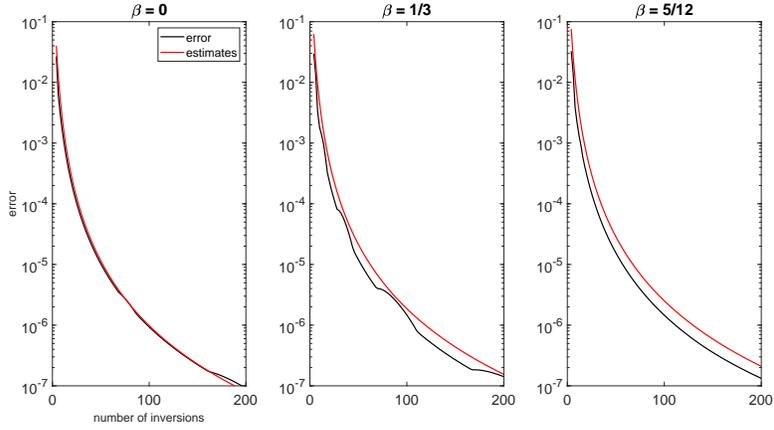}
\end{center}
\caption{Error and error estimate (\ref{bound1}) for $\beta=0$, $\beta=\frac{1}{3}$, $\beta=\frac{5}{12}$, from left to right.}
\label{Figura3}
\end{figure}

In order to test the behavior of the method, we consider a diagonal test matrix with a very large spectrum. In particular, we define
\begin{equation*}
\mathcal{L} = {\rm diag} \left(1, 1+\rho_1 e^{i \beta \pi}, 1+\rho_1 e^{-i \beta \pi}, \ldots, 1+\rho_N e^{i \beta \pi}, 1+\rho_N e^{-i \beta \pi} \right),
\end{equation*}
where 
\begin{equation*}
\rho_i:=10^{x_i}, \quad x =\left(0,0.1,0.2, \ldots, 16\right).
\end{equation*}
The matrix is clearly normal so that $F(\mathcal{L})$ is the convex hull of the spectrum, that is, the triangle with vertex at $1$, $1+10^{16} e^{i \beta \pi}$, $1+10^{16} e^{-i \beta \pi}$. 
In Figure \ref{Figura3} we plot the error toghether with the estimate (\ref{bound1}) for some values of $\beta$.
In Table \ref{tabella1} we also show the corresponding values assumed by $\tau$ and $\hat{\rho}$ relative to the rightmost plot of Figure \ref{Figura3} ($\beta=5/12$).
\begin{table}[tbp]
\begin{equation*}
\begin{array}{cccccccc}
\toprule
n & 10 & 25 & 40 & 55 & 70 & 85 & 100 \\
\midrule
\tau & 1.0E02 & 1.1E03  & 4.2E03 & 1.1E04 & 2.3E04 & 4.1E04 & 6.8E04 \\
\hat{\rho} & 4.1E05 & 2.2E08 & 6.0E09 & 5.7E10 & 3.1E11 & 1.3E12 & 4.0E12 \\
\bottomrule
\end{array}
\end{equation*}%
\caption{Values of $\tau$ and $\hat{\rho}$ for $\beta=\frac{5}{12}$.}
\label{tabella1}
\end{table}

\section{The case of bounded operators} \label{Section5}

Consider the case of bounded sectorial operators $\mathcal{L}_N$ with numerical range contained in 
\begin{equation*}
\Sigma_{\beta,1,\rho_N} =\left\lbrace z \in \mathbb{C} \mid z=1+\rho e^{i \theta \pi}, | \theta | \leq \beta,\, \beta < \frac{1}{2}, \,  0 \leq \rho \leq \rho_N \right\rbrace.
\end{equation*}
Before starting, we need to remember that the local maximum $\rho^{\star}$ of $\Phi^{(1)}\left(\tau, 1+\rho e^{i \beta \pi}\right)$ grows approximatively like $n^8/(\ln n)^4$ (see (\ref{rho_star_exp})). Therefore, there exists $\bar{n}$ such that $\rho^{\star}>\rho_N$ for $n>\bar{n}$. As consequence, for $\rho^{\star} \leq \rho_N$, the estimate (\ref{bound1}) is still valid, because $\Phi^{(1)}\left(\tau, 1+\rho^{\star}e^{i \beta \pi}\right) \geq \Phi^{(1)}\left(\tau, 1+\rho_N e^{i \beta \pi}\right)$, and hence we have to solve (\ref{condizione tau}) to approximate the solution of (\ref{condizionePHI}). On the other side, for $\rho^{\star} > \rho_N$ the bound can be improved, because the $\rho^{\star}$ falls outside $[0,\rho_N]$.
Similarly to the analysis given in \cite{AcetoNovati}, for $\rho^{\star} > \rho_N$ the optimal value for $\tau$ can be approximated by solving
\begin{equation} \label{condizione tau bounded}
g^{(1)}(\tau, \rho_N) = g^{(2)} (\tau,\rho_0).
\end{equation}

\begin{proposition}
For $1 \ll \tau \ll \rho_N$, the solution of (\ref{condizione tau bounded}) is approximated by
\begin{equation} \label{tau_hat}
\hat{\tau}= \left( -\frac{\rho_N^{\frac{1}{4}}}{8\sqrt{2}C_{\beta} n} \ln \left(\frac{\sqrt{\rho_N}}{D_{\beta}}\right)  + \sqrt{\left(\frac{\rho_N^{\frac{1}{4}}}{8 \sqrt{2}C_{\beta} n} \ln \left(\frac{\sqrt{\rho_N}}{D_{\beta}}\right) \right)^2+\frac{G_{\beta}}{C_{\beta}} \rho_N^{\frac{1}{4}}}\right)^4.
\end{equation}
\end{proposition}
\begin{proof}
Using relations (\ref{g1}) and (\ref{C}), equation (\ref{condizione tau bounded}) becomes
\begin{equation} \label{D}
\rho_N^{-\frac{1}{2}} \left[ 1+C_{\beta} \sqrt{2} \left(\frac{\tau}{\rho_N}\right)^{\frac{1}{4}} \right]^{-2n}= \frac{1}{D_{\beta}} \left(1+\sqrt{2}G_{\beta}\tau^{-\frac{1}{4}}\right)^{-2n}.
\end{equation}
Using the approximations
\begin{equation} \label{approx C_N}
1+C_{\beta} \sqrt{2} \left(\frac{\tau}{\rho_N}\right)^{\frac{1}{4}} \cong \exp \left(C_{\beta} \sqrt{2} \left(\frac{\tau}{\rho_N}\right)^{\frac{1}{4}}\right),
\end{equation}
and
\begin{equation*}
1+\sqrt{2}G_{\beta}\tau^{-\frac{1}{4}} \cong \exp \left(\sqrt{2}G_{\beta}\tau^{-\frac{1}{4}}\right),
\end{equation*}
we rewrite the equation (\ref{D}) as
\begin{equation*}
\rho_N^{\frac{1}{4n}} \exp \left( C_{\beta} \sqrt{2} \left(\frac{\tau}{\rho_N}\right)^{\frac{1}{4}} \right) = D_{\beta}^{\frac{1}{2n}} \exp \left( G_{\beta}\sqrt{2}\tau^{-\frac{1}{4}} \right).
\end{equation*}
Therefore, 
\begin{equation*}
\tau^{\frac{1}{2}}+\frac{\rho_N^{\frac{1}{4}}}{4C_{\beta} \sqrt{2} n} \ln \left( \frac{\sqrt{\rho_N}}{D_{\beta}} \right) \tau^{\frac{1}{4}} -\frac{G_{\beta}}{C_{\beta}}\rho_N^{\frac{1}{4}}=0.
\end{equation*}
By solving this equation and taking the positive solution we obtain the result.
\end{proof}

We observe that by (\ref{tau_hat}), for $n \rightarrow +\infty$ we have
\begin{align*}
\left(\frac{\hat{\tau}}{\rho_N}\right)^{\frac{1}{4}} &= - \frac{1}{4  \sqrt{2}C_{\beta}n} \ln \left(\frac{\sqrt{\rho_N}}{D_{\beta}}\right)+\sqrt{\left(\frac{1}{4 \sqrt{2}C_{\beta}n}\ln \left(\frac{\sqrt{\rho_N}}{D_{\beta}}\right) \right)^2+ \frac{G_{\beta}}{C_{\beta}}\rho_N^{-\frac{1}{4}}} \\
&= -\frac{1}{4 \sqrt{2}C_{\beta}n}\ln \left(\frac{\sqrt{\rho_N}}{D_{\beta}}\right)+\sqrt{\frac{G_{\beta}}{C_{\beta}}} \rho_N^{-\frac{1}{8}}  + \mathcal{O} \left(\frac{1}{n^2} \right).
\end{align*}
Using (\ref{Ecerchiato}), (\ref{k_1 approx}), (\ref{approx C_N}) in order and the above result, for $n >\bar{n}$ we obtain
\begin{align}
|E_n(\mathcal{L}_N)| &\cong 4 \rho_N^{-\frac{1}{2}} \left(1+\sqrt{2} C_{\beta} \left( \frac{\hat{\tau}}{\rho_N} \right)^{\frac{1}{4}} \right)^{-2n} \notag \\
&\cong 4 \rho_N^{-\frac{1}{2}} \exp\left(-2 \sqrt{2} C_{\beta} n \left( \frac{\hat{\tau}}{\rho_N} \right)^{\frac{1}{4}}\right) \notag \\
& \cong 4 \rho_N^{-\frac{1}{4}}D_{\beta}^{-\frac{1}{2}} \exp\left(-2 \sqrt{2}\sqrt{G_{\beta}C_{\beta}} n \rho_N^{-\frac{1}{8}}\right). \label{errore_bound}
\end{align}
Note that $\sqrt{G_{\beta}C_{\beta}}=1$ for $\beta=0$.

In order to derive an estimate of $\bar{n}$, we impose $\hat{\rho} = \rho_N$, where $\hat{\rho}$ is as in (\ref{rho_star_exp}). We obtain the equation
\begin{equation*}
\exp \left( W(2eC_{\beta}(n-1)n) \right) =e \rho_N^{\frac{1}{4}},
\end{equation*}
and therefore
\begin{equation*}
W(2eC_{\beta}(n-1)n) = \ln \left(e \rho_N^{\frac{1}{4}} \right).
\end{equation*}
Since $W(z_1) =z_2$ if and only if $z_1=z_2e^{z_2}$, it follows that
\begin{align*}
\bar{n} &\cong \frac{\rho_N^{\frac{1}{8}}}{\sqrt{2C_{\beta}}} \left(\ln \left(e \rho_N^{\frac{1}{4}} \right) \right)^{\frac{1}{2}} \\
&\cong \frac{\rho_N^{\frac{1}{8}} \left(\ln \rho_N \right)^{\frac{1}{2}}}{2 \sqrt{2 C_{\beta}}}.
\end{align*}

\subsection{Numerical experiments}

\begin{figure}
\begin{center}
\includegraphics[scale=0.3]{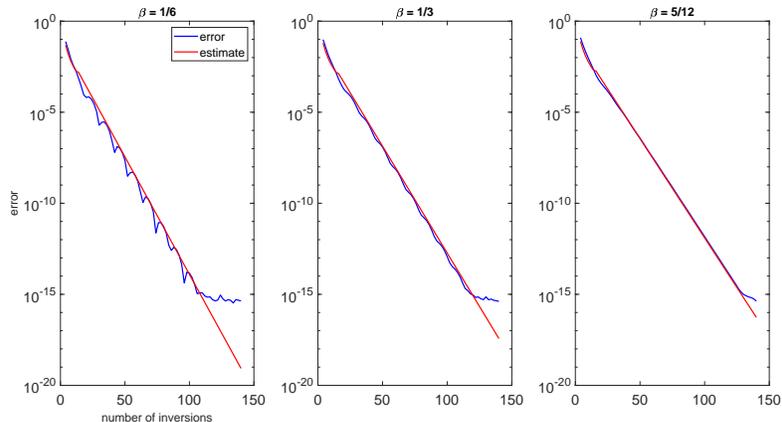}
\end{center}
\caption{Error and error estimate (\ref{errore_bound}) for $\beta=\frac{1}{6}$, $\beta=\frac{1}{3}$ and $\beta=\frac{5}{12}$, from left to right.}
\label{Figura7}
\end{figure}

In order to test the method for sectorial bounded operators, we first consider the same operator of Section \ref{Section4.4} but with
\begin{equation*}
x=\left(0,0.1,0.2, \ldots, 4\right),
\end{equation*}
so that $\rho_N=10^4$. In Figure \ref{Figura7} the error and the estimate (\ref{errore_bound}) are plotted for different values of $\beta$.

\begin{figure}
\begin{center}
\includegraphics[scale=0.3]{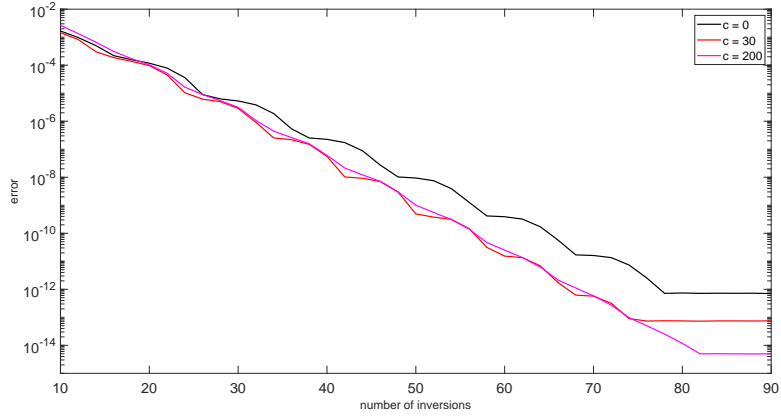}
\end{center}
\caption{Error for the discretization of the operator defined in (\ref{operatore_derivara}), with $c=0$, $c=30$ and $c=200$.}
\label{Figura5}
\end{figure}

As a more realistic example, we also consider the discretization using central differences of the operator
\begin{equation} \label{operatore_derivara}
\mathcal{L} u=-u{''} + c u{'}, \quad c \geq 0, 
\end{equation}
on $[0,1]$ with Dirichlet boundary conditions. We have taken $N=200$ equally spaced interior points. By moving the constant $c$ we change the angle $\beta \pi$ of the sector containing $F(\mathcal{L}_N)$, where $\mathcal{L}_N$ is the discretization matrix.
In Figure \ref{Figura5} we plot the errors for $c=0$ ($\beta =0$), $c=30$ ($\beta=0.44$) and $c=200$ ($\beta=0.49$). It is interesting to observe that the method is a bit faster for $c \gg 0$. This is due to the position of the eigenvalues of smallest modulus that move away from $0$ for growing $c$. With large $c$ we also notice an improvement of the attainable accuracy and the reason lies in the conditioning of $\mathcal{L}_N$ that reduces incresing $c$.

\section{Conclusions}

We have studied an unexplored Gaussian approach for the computation of the inverse square root of regularly accretive operators.
The method exhibits a very fast initial convergence and its rate is almost independent of the angle of the sector containing the numerical range. We have derived sharp error estimates that can be used for an a priori selection of the number of quadrature points necessary to achieve a prescribed accuracy.

\begin{figure}
\begin{center}
\includegraphics[scale=0.3]{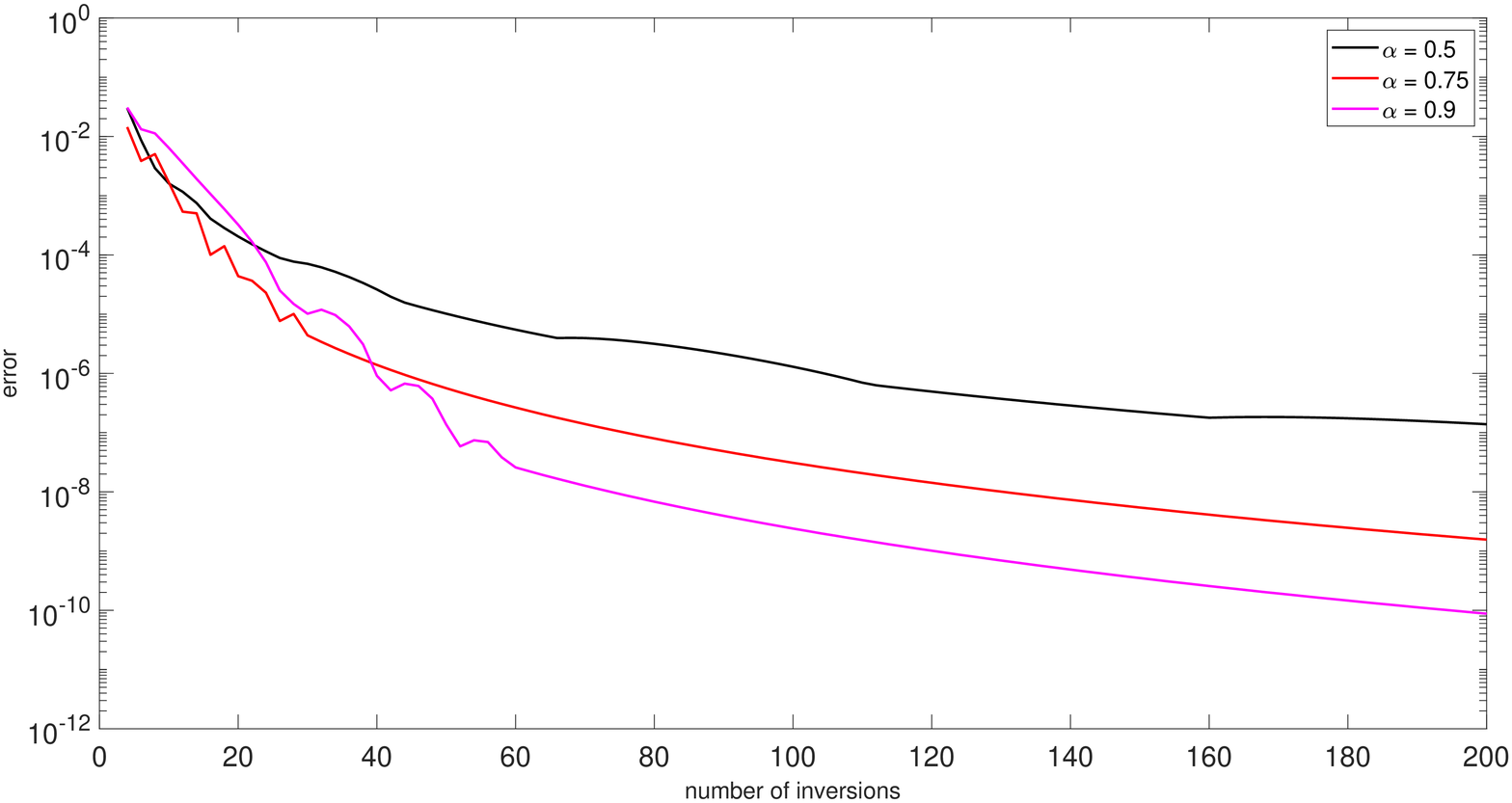}
\end{center}
\caption{Error for the operator $\mathcal{L}^{-\alpha}$ with $\beta=\frac{1}{6}$, for $\alpha=0.5$, $\alpha=0.75$ and $\alpha=0.9$.}
\label{Figura8}
\end{figure}

While all the analysis is restricted to the square root, we remark that the method works fine also more generally for $\mathcal{L}^{-\alpha}$, $0<\alpha<1$, and in particular for $0.5<\alpha<1$. 
Working with the example of Section \ref{Section4.4}, in Figure \ref{Figura8} we compare the behavior of the method for $\alpha=0.75$ and $\alpha=0.9$ with respect to the case $\alpha=0.5$, using the value of $\tau$ derived in Section \ref{Section4.2}.
The reason for which we have not studied the general case of $\alpha \neq 0.5$ lies in the regularity of the integrand function arising from the changes of variable we have used at the beginning. Indeed, following the same approach, for $0<\alpha<1$ we obtain
\begin{equation*}
\lambda^{-\alpha} = \frac{\sin (\alpha \pi)}{\pi} \tau^{1-\alpha} \left[ \frac{1}{\alpha} \int_0^1 \frac{1}{\tau +\lambda t^{\frac{1}{\alpha}}} dt+\frac{1}{1-\alpha} \int_0^1 \frac{1}{\tau t^{\frac{1}{1-\alpha}}+\lambda}dt \right].
\end{equation*}
Consequently, the integrand functions are simultaneously analytic only for $\alpha=0.5$. Taking for istance $\alpha= \frac{1}{k}$, $k \in \mathbb{N}$, the second integrand function is not analytic at $0$, so that all the error analysis becomes extremely complicated (see e.g. \cite{Barrett},\cite[Section 4]{Trefethen}).

\appendix

\section{Approximation of $\beta^{\star}$}

Let $\Psi_{s_0}\in \mathcal{E}$, be the ellipse passing through the point $\frac{2}{\sqrt{\tau}}i-1$ (cf. (\ref{22_bis})).
The value $\beta^{\star}$, such that the function $\Phi^{(2)} \left(1+\rho e^{-i \beta \pi} \right)$ possesses a local maximum for $\beta>\beta^{\star}$, is the one for which $\Psi_{s_0}$ is also tangent to the curve $\chi^{(2)}\left(\Gamma_{\beta^{\star}}^- \right)$ at $\frac{2}{\sqrt{\tau}}i-1$ (see Figure \ref{Figura6}). In order to compute $\beta^{\star}$, we consider the tangents at $\frac{2}{\sqrt{\tau}}i-1$ to the ellipse and to the curve, and impose them to have the same slope. Before starting, we need to derive the semi-width $\gamma$ and the semi-height $\delta$ of the ellipse. By geometrical evidence, we have that
\begin{equation} \label{delta}
\delta = \Im \left\lbrace \frac{1}{2} \left( s_0 e^{i \frac{\pi}{2}}+\frac{1}{s_0} e^{-i \frac{\pi}{2}} \right) \right\rbrace = \frac{1}{2} \left( s_0-\frac{1}{s_0} \right),
\end{equation}
and
\begin{equation*}
\gamma^2 = \delta^2+1.
\end{equation*}
At this point, we remind that the slope of the tangent at $\frac{2}{\sqrt{\tau}}i-1$ to the ellipse is
\begin{align}
m &= -\frac{\delta^2}{\gamma^2}\frac{\Re\left(\frac{2}{\sqrt{\tau}}i-1 \right)}{\Im\left(\frac{2}{\sqrt{\tau}}i-1 \right)} \notag \\
&=\frac{\delta^2}{\delta^2 +1} \left( \frac{\sqrt{\tau}}{2}\right) \notag \\
&\sim  \frac{\sqrt{\tau}}{2+\sqrt{\tau}} \sim 1, \quad \tau \rightarrow + \infty, \label{slope}
\end{align} 
where we have used (\ref{delta}) and 
\begin{equation*}
s_0-\frac{1}{s_0} \sim \frac{4}{\sqrt{2}} \tau^{-\frac{1}{4}},
\end{equation*} 
that comes from (\ref{22_bis}).
Now, it is not difficult to show that the angle between the tangent to the curve $\chi^{(2)}\left(\Gamma_{\beta}^- \right)$ at $\frac{2}{\sqrt{\tau}}i-1$ and the line $\Im (z) = \frac{2}{\sqrt{\tau}}$ is given by $\frac{1}{2} \left(\pi-2\beta \pi\right)$.
Hence, in order to find an approximation of $\beta^{\star}$, by (\ref{slope}) we impose the condition 
\begin{equation*}
\tan \left[ \frac{\pi}{2} (1-2\beta) \right] = 1,
\end{equation*}
that leads to $\beta^{\star} \sim \frac{1}{4}$ as $\tau \rightarrow +\infty$.

\section*{Acknowledgements}

This work was partially supported by GNCS-INdAM, FRA-University of Trieste and CINECA under HPC-TRES program award number 2019-04. The authors are members of the INdAM research group GNCS.

\end{document}